\documentclass[MSNbibl,number,citesort,dvips]{arxbj}
\usepackage{mathbh}

\aid{0}
\volume{20}
\issue{2}
\pubyear{2014}
\firstpage{979}
\lastpage{989}
\doi{10.3150/13-BEJ512} 

\makeatletter
\def\cal{\mathcal}
\newtheorem{theorem}{Theorem}
\newtheorem{corollary}[theorem]{Corollary}
\newtheorem{lemma}[theorem]{Lemma}
\def\L{{\cal L}}
\def\p{\mathbb{P}}
\def\R{\mathbb{R}}
\def\a{\alpha}
\def\E{\mathbb{E}}
\def\ve{\varepsilon}
\def\1{\mathbh{1}}
\def\d{d_{\mathrm{TV}}}
\makeatother

\begin{document}
\begin{frontmatter}

\title{Lower bounds to the accuracy of inference on heavy tails}
\runtitle{Lower bounds}

\begin{aug}
\author{\fnms{S.Y.}~\snm{Novak}\corref{}\ead[label=e1]{S.Novak@mdx.ac.uk}}
\runauthor{S.Y. Novak} 
\address{Middlesex University, The Burroughs, London NW44BT, UK.
\printead{e1}}
\end{aug}

\received{\smonth{6} \syear{2012}}
\revised{\smonth{1} \syear{2013}}

%
\begin{abstract}
The paper suggests a simple method of deriving minimax lower bounds to
the accuracy of statistical inference on heavy tails.
A well-known result by Hall and Welsh (\textit{Ann. Statist.} \textbf{12} (1984) 1079--1084) states that if $ \hat
\a_n $ is an estimator of the tail index $ \a_P $ and $ \{z_n\} $
is a sequence of positive numbers such that $ \sup_{P\in{\cal D}_r}
\p(|\hat\a_n - \a_P| \ge z_n) \to0$, where $ {\cal D}_r $ is
a certain class of heavy-tailed distributions, then $ z_n\gg n^{-r}$.
The paper presents a non-asymptotic lower bound to the probabilities
$ \p(|\hat\a_n-\a_P| \ge z_n)$. We also establish non-uniform lower
bounds to the accuracy of tail constant and extreme quantiles
estimation. The results reveal that normalising sequences of robust
estimators should depend in a specific way on the tail index and the
tail constant.
\end{abstract}

%

\end{frontmatter}

\section{Introduction}

A growing number of publications is devoted to the problem of statistical
inference on heavy-tailed distributions. Such distributions naturally appear
in finance, meteorology, hydrology, teletraffic engineering, etc. \cite{EKM,R97}.
In particular, it is widely accepted that frequent financial data (e.g.,
daily and hourly log-returns of share prices, stock indexes and
currency exchange rates) often exhibits heavy tails \cite
{FR,EKM,M63,N09}, while less frequent financial data is typically light-tailed.
The heaviness of a tail of the distribution appears to be responsible for
extreme movements of stock indexes and share prices.
The tail index indicates how heavy the tail is; extreme quantiles are
used as measures of financial risk \cite{EKM,N09}. The need to evaluate
the tail index and extreme quantiles stimulated research on methods of
statistical inference on heavy-tailed data.

The distribution of a random variable (r.v.) $ X $ is said to have a
\textit{heavy right tail} if
%
\begin{equation}
\label{T1} \p(X\ge x) = L(x)x^{-\a} \qquad (\a>0),
\end{equation}
where the (unknown) function $ L $ is slowly varying at infinity:
\[
\lim_{x\to\infty} L(xt)/L(x) = 1 \qquad(\forall t>0) .
\]
We denote by $ {\cal H} $ the non-parametric class of distributions
obeying (\ref{T1}).\vadjust{\goodbreak}

The tail index $ \a$ is the main characteristic describing the tail
of a
distribution. If $ L(x) = c+\mathrm{o}(1),$ then $ c $ is called the tail constant.

Let $ F(\cdot)=\p(X < \cdot) $ denote the distribution function (d.f.).
Obviously, the tail index is a functional of the distribution function:
%
\begin{equation}
\label{T3} \a_{F}\equiv\a_{P} = -\lim_{x\to\infty}
\frac{\ln\p(X\ge
x)}{\ln
x} .
\end{equation}
%
If $ L(x) $ tends to a constant (say, $c_{F}$) as $ x \to\infty,$
then the tail constant is also a functional of $ F$:
\[
c_{F}\equiv c_{P} = \lim_{x\to\infty}
x^{\a_F} \p(X\ge x) .
\]

The statistical inference on a heavy-tailed distribution is
straightforward if the class of unknown distributions is assumed to be
a regular parametric family. The drawback of the parametric approach is
that one usually cannot reliably check whether the unknown distribution
belongs to a chosen parametric family.

A lot of attention during the past three decades has been given to
the problem of reliable inference on heavy tails without parametric assumptions.
The advantage of the non-parametric approach is that a class of
unknown distributions, $ {\cal P},$ is so large that the problem of
testing the hypothesis that the unknown distribution belongs to $ {\cal
P} $ does not arise.
The disadvantage of the non-parametric approach is that virtually no
question concerning inference on heavy tails can be given a simple answer.
In particular, the problem of establishing a lower bound to the
accuracy of
tail index estimation remained open for decades.

A lower bound to the accuracy of statistical inference sets a
benchmark against which the accuracy of any particular estimator can be
compared.
When looking for an estimator $ \hat a_n $ of a quantity of interest,
$ a_P,$ where $ P\in{\cal P} $ is the unknown distribution, $ {\cal P}
$ is the class of distributions and $ a_P $ is a functional of $ P,$
one often would like to choose an estimator that minimises a loss
function uniformly in $ {\cal P} $ (e.g., $ \sup_{P\in{\cal P}} \E_P
\ell(|\hat a_n-a_P|),$ where $\ell$ is a particular loss function).
A~lower bound to $ \sup_{P\in{\cal P}} \E_P \ell(|\hat a_n-a_P|) $
follows if one can establish a lower bound to
\[
\sup_{P\in{\cal P}} \p\bigl(|\hat a_n-a_P|\ge u\bigr)\qquad
(u>0).
\]

The first step towards establishing a lower bound to the accuracy of tail
index estimation was made by Hall and Welsh \cite{HW84}, who proved the
following result. Note that the class $ {\cal H} $ of heavy-tailed
distributions is too ``rich'' for meaningful inference, and one usually
deals with a subclass of $ {\cal H},$ imposing certain restrictions on
the asymptotics of $ L(\cdot)$.
Hall and Welsh dealt with the class $ {\cal D}_{b,A} \equiv{\cal
D}_{b,A} (\a_0,c_0,\ve) $ of distributions on $ (0;\infty) $ with densities
%
\begin{equation}
\label{Tdens} f(x) = c\a x^{-\a-1}\bigl(1+u(x)\bigr),
\end{equation}
where $ \sup_{x>0}|u(x)|x^{b\a}\le A,$ $ |\a-\a_0|\le\ve,$ $ |c-c_0|
\le\ve$.
Note that the range of possible values of the tail index is restricted
to interval $ [\a_0-\ve;\a_0+\ve]$. Let
\[
\hat\a_n \equiv\hat\a_n(X_1,\ldots,X_n)
\]
be an arbitrary tail index estimator, where $ X_1,\ldots,X_n $ are
independent and identically distributed (i.i.d.) random variables, and
let $ \{z_n\} $ be a sequence of positive numbers. If
%
\begin{equation}
\label{HW} \lim_{n\to\infty} \sup_{F\in{\cal D}_{b,A}}
\p_F\bigl(|\hat\a_n-\a _F| \ge z_n\bigr)=0\qquad
(\forall A>0),
\end{equation}
then
\[
z_n \gg n^{-b/(2b+1)}\qquad (n\to\infty)
\]
(to be precise, Hall and Welsh \cite{HW84} dealt with the random variables
$ Y_i=1/X_i,$ where $ X_i $ are distributed according to (\ref{Tdens})).

Beirlant et al. \cite{BBW06} have a similar result for a larger class $ \cal
P $ of distributions but require the estimators are uniformly
consistent in $ \cal P$.
Pfanzagl \cite{Pf} has established a lower bound in terms of a modulus
of continuity related to the total variation distance $\d$.
Let $ {\cal D}_b^+ $ be the class of distributions with densities
(\ref{Tdens}) such that $ \sup_{x>0} |u(x)|x^{\a b}<\infty, \a>0$,
and set
\[
s_n(\ve,P_0) = \sup_{P\in{\cal P}_{n,\ve}} |
\a_P-\a_{P_0}|,
\]
where $ \a_P $ is the tail index of distribution $ P $ and
$ {\cal P}_{n,\ve} = \{P\in{\cal D}_b^+\dvt \d(P_0^n;P^n)\le\ve\} $
is a
neighborhood of $ P_0\in{\cal D}_b^+ .$
Pfanzagl has showed that neither estimator can converge to $ \a$ uniformly
in $ {\cal P}_{n,\ve} $ with the rate better than $ s_n(\ve,P_0),$ and
\[
\inf_{0<\ve<1}\ve^{-2b/(1+2b)}\liminf_{n\to\infty}
n^{b/(1+2b)} s_n(\ve ,P_0)>0 .
\]
Donoho and Liu \cite{DL} present a lower bound
to the accuracy of tail index estimation in terms of a modulus of
continuity $ \Delta_A(n,\ve)$. However, they do not calculate $
\Delta
_A(n,\ve)$. The claim that a particular heavy-tailed distribution is
stochastically dominant over all heavy-tailed distributions with the
same tail index appears without proof.
Assuming that the range of possible values of the tail index is
restricted to an interval of fixed length, Drees \cite{Drees2001}
derives the asymptotic minimax risk for affine estimators of the tail
index and indicates an approach to numerical computation of the
asymptotic minimax risk for non-affine ones.

The paper presents a simple method of deriving minimax lower bounds to
the accuracy of non-parametric inference on heavy-tailed distributions.
The results are non-asymptotic, the constants in the bounds are shown
explicitly, the range of possible values of the tail index is not
restricted to an interval of fixed length. The information functional
seems to be found for the first time, as well as the lower bound to the
accuracy of extreme quantiles estimation.

The results indicate that the traditional minimax approach may require
revising. The classical approach suggests looking for an estimator $
\hat a_n $ that minimises, say,
\[
\sup_{P\in{\cal P}} \E_P |\hat a_n -
a_P|
\]
(cf. \cite{Hu97,IH81,Tsy}), while our results suggest looking for an
estimator $ \hat a^*_n $ that minimises
\[
\sup_{P\in{\cal P}} g_P \E_P\bigl|\hat
a_n^*-a_P\bigr|,
\]
where $ g_P $ is the ``information functional'' (an analogue of Fisher's
information).
Theorems \ref{LP-1}--\ref{LP-3} reveal the information functionals and
indicate that the normalising sequence of a robust estimator should
depend in a specific way on the characteristics of the unknown
distribution.\vspace*{-2pt}

\section{Results}\vspace*{-2pt}

In the sequel, we deal with the non-parametric class
%
\begin{equation}
\label{T4} {\cal H}(b) = \Bigl\{ P\in{\cal H}\dvt \sup_{x>K_*(P)} \bigl|
c_{F}^{-1} x^{\a_{F}}P(X\ge x) - 1 \bigr| x^{b\a_{F}}
< \infty\Bigr\} 
\end{equation}
of distributions on $ (0;\infty)$, where $ b>0 $ and $ K_*(P) $ is the
left end-point of the distribution. If $ \L(X)\in{\cal H}(b),$ then
\[
\p(X\ge x) = c_{F} x^{-\a_{F}} \bigl(1+\mathrm{O}\bigl(x^{-b\a_{F}}
\bigr) \bigr)\qquad (x\to\infty) .
\]
The class $ {\cal H}(b) $ is larger than $ {\cal D}_b^+$; the range of
possible values of the tail index is not restricted to an interval of
fixed length. Below, given a distribution function (d.f.) $ F_i$, we
put\looseness=-1
\[
a_{{F_i}}=1/\a_{{F_i}} ,\qquad  r=b/(1+2b) ,
\]\looseness=0
$\E_i $ means the mathematical expectation with respect to $ F_i $ and
$ \p_i $ is the corresponding distribution. We set $ K\equiv K_{\a,b,c}
= \a^{-2r}c^{-\a r}\mathrm{e}^{-1}  (c^{\a b} \wedge \mathrm{e}^{-2b} ) $.

\begin{theorem} \label{LP-1} 
For any $ \a>0,$ $ c>0,$ any tail index estimator $ \hat\a_n $ and any
estimator $ \hat a_n $ of index $ a=1/\a$ there exist d.f.s $ F_0,
F_1\in{\cal H}(b) $ such that $ \a_{{F_0}}=\a, c_{{F_0}}=c^{-\a
},$ and
%
\begin{eqnarray}
\label{LP1} \max_{i\in\{0;1\}} \p_i\bigl(|\hat
\a_n/\a_{{F_i}}-1| \a_{{F_i}}^{r/b}
c_{{F_i}}^r n^r \ge v/2\bigr) &\ge&
\bigl(1-v^{1/r}/8n\bigr)^{2n}/4 ,
\\
\label{LP} \max_{i\in\{0;1\}} \p_i\bigl(|\hat
a_n/a_{{F_i}}-1| a_{{F_i}}^{-r/b}c_{{F_i}}^r
n^r \ge v/2\bigr)& \ge&\bigl(1-v^{1/r}/8n\bigr)^{2n}
/4
\end{eqnarray}
as $ n>4\max\{\a^2c^{-2\a b};c^{2\a}\a^{-2/b}\} $ and $ v\in[0;K n^r]$.
\end{theorem}

Note that if $ \max_{i\in\{0;1\}} \p_i(|\hat\a_n/\a_{{F_i}}-1|
\ge
z_n) \to0 $ as $ n\to\infty,$ then for any $ C>0 $ we have $ z_n\ge
Cn^{-r} $
for all large enough $ n$, yielding $ z_n\gg n^{-r} $. Thus, the
Hall--Welsh result follows from~(\ref{LP1}).

Theorem \ref{LP-1} shows that the natural normalising sequence for $
\hat\a_n/\a_{F}-1 $ is $ n^{-r}\a_{F}^{-r/b}c_{F}^{-r} .$
The information functional $ g_F = \a_{F}^{r/b} c_{{F}}^r $ plays
here the same role as Fisher's information function in the Fr\'
echet--Rao--Cram\'er inequality.

Theorem \ref{LP-1} yields also minimax lower bounds to the moments of
$ |\hat\a_n/\a_{{F_i}}-1|.$ In particular, there holds\vadjust{\goodbreak}

\begin{corollary} \label{LP-C}
For any $ \a>0,$ $ c>0 $ there exist distribution functions $
F_0,F_1\in
{\cal H}(b) $ such that $ \a_{{F_0}} = \a, c_{{F_0}} = c^{-\a},$ and
for any tail index estimator $ \hat\a_n $ 
%
\begin{equation}
\label{LPC} \max_{i\in\{0;1\}} \a_{{F_i}}^{r/b}
c_{{F_i}}^r \E_{{F_i}} |\hat\a_n/
\a_{{F_i}}-1|n^{r} \ge4^{r}r\Gamma(r)/8+\mathrm{o}(1) .
\end{equation}
The result holds if $ \a_{{F_i}}^{r/b}c_{{F_i}}^r \E_{{F_i}} |\hat
\a
_n/\a_{F_i}-1| $ in (\ref{LPC}) is replaced with $ a_{{F_i}}^{-r/b}
c_{{F_i}}^r \E_{{F_i}} |\hat a_n/a_{{F_i}}-1|.$
\end{corollary}

Let $ {\cal H}_n(b)\subset{\cal H}(b) $ be a class of d.f.s such that
$ \inf_{F\in{\cal H}_n(b)} K_{\a_{F},b,c_{F}} n^r \to\infty$ as $
n\to\infty$. Then for any estimator $ \hat\a_n $
{\renewcommand{\theequation}{$\protect\ref{LPC}^*$}
\begin{equation}\label{LPCast}
\sup_{F\in{\cal H}_n(b)} \a_{F}^{r/b}
c_{{F}}^r \E_{{F}} |\hat\a_n/
\a_{F}-1|n^{r} \ge4^{r}r\Gamma(r)/8+\mathrm{o}(1).
\end{equation}}
\hspace*{-2pt}A lower bound to $ \E_{{F}} |\hat\a_n/\a_{F}-1| $ seems to be
established for the first time.

The presence of the information functional makes the bound
non-uniform. Note that a uniform lower bound would be meaningless: as
the range of possible values of $ \a_{F} $ is not restricted to an
interval of fixed length, it follows from (\ref{LPCast}) that
\[
\sup_{F\in{\cal H}_n(b)} \E_F |\hat\a_n/
\a_{F}-1| \to\infty\qquad (n\to \infty).
\]
More generally, $ \sup_{F\in{\cal H}_n(b)} \tilde g_F\E_F |\hat\a
_n/\a
_{F}-1| $ may tend to $ \infty$ as $ n\to\infty$ if $ \tilde
g_F/g_F\ne \mathrm{const}$.

Let $ \hat c_n $ be an arbitrary tail constant estimator.
The next theorem presents a lower bound to the probabilities $ \p_F(|\hat
c_n-c_{F}| \ge x)$.

\begin{theorem} \label{LP-2}
Let $ \hat c_n $ be an arbitrary tail constant estimator. For any $ \a
\ge n^{-r/2} $ and $ c>0 $ there exist distribution functions $
F_0,F_1\in{\cal H}(b) $ such that $ \a_{{F_0}} = \a, c_{{F_0}} =
c^{-\a},$ and for all large enough $ n$, as $ v\in[0;\a^{-2}c^{-\a
}\ln n]$,
%
\setcounter{equation}{8}
\begin{equation}
\label{LP4} \max_{i\in\{0;1\}} \p_i \bigl( |\hat
c_n/c_{{F_i}} - 1|\a_{{F_i}}^{r/b}
c_{{F_i}}^r \ge rv^rn^{-r}\ln(n/\ln
n)t_n/2b \bigr) \ge(1 - v/8n)^{2n}/4 ,
\end{equation}
where $ t_n=\exp(-r(1 - r)n^{-r/2}(\ln(n/\ln n))^{r+1}/b).$
\end{theorem}

Similarly to (\ref{LPC}) Theorem \ref{LP-2} yields lower bounds to the
moments of $ |\hat c_n/c_{{F_i}}-1|.$ In particular, (\ref{LP4}) entails
{\renewcommand{\theequation}{$\protect\ref{LP4}^*$}
\begin{equation}\label{LP4ast}
\max_{i\in\{0;1\}} \a_{{F_i}}^{r/b}
c_{{F_i}}^r \E_{{F_i}} |\hat c_n/c_{F_i}-1|
\ge(\ln n)n^{-r}r^24^{r-1}\Gamma(r)/\bigl(2b+\mathrm{o}(1)
\bigr) . 
\end{equation}}

According to Hall and Welsh \cite{HW84},
\[
z_n \gg(\ln n) n^{-b/(2b+1)}
\]
if $ \lim_{n\to\infty} \sup_{F\in{\cal D}_{b,A}} \p_F(|\hat c_n-c_F|
\ge z_n)=0\
(\forall A>0) $. This fact can be obtained as a consequence to Theorem
\ref{LP-2}: if $ \max_{i\in\{0;1\}} \p_i(|\hat c_n-c_{{F_i}}| \ge z_n)
\to0 $ as $ n\to\infty,$ then for any $ C>0 $ we have $ z_n \ge
Cn^{-r}\ln n $ for all large enough $ n$, hence $ z_n \gg n^{-r}\ln
n$.

We now present a lower bound to the accuracy of estimating extreme upper
quantiles. We call an upper quantile of level $ q $ ``extreme'' if $
q\equiv q_n $ tends to 0 as $ n $ grows.
In financial applications (see, e.g., \cite{EKM,N09}), an upper
quantile of the level as high as 0.05 can be considered extreme as the
empirical quantile estimator appears unreliable. Of course, there is an
infinite variety of possible rates of decay of $ q_n.$ Theorem \ref
{LP-3} presents lower bounds in the case $ q_n = s n^{-1/(1+2b)},$
where $ s $ is bounded away from $ 0 $ and $ \infty$.

Set $ \bar F = 1-F.$ We denote the upper quantile of level $ q_n $ by
\[
x_{F,n} = \bar F^{-1}(q_n) .
\]
Let $ \hat x_n $ be an arbitrary estimator of $ x_{F,n}$. Denote
$ w_{F_i} \equiv w_{F_i}(\a_{F_i},c_{F_i},b,s,u)
=\break   |\ln( u\a_{{F_i}}^{2r} c_{{F_i}}^{2br} /s^b )|$.

\begin{theorem} \label{LP-3}
For any $ \a>0,$ $ c>0 $ there exist distribution functions
$ F_0,F_1 \in{\cal H}(b) $ such that $ \a_{{F_0}} = \a,
c_{{F_0}}=c^{-\a},$ and for all large enough $ n $ and $ u\in( s^b\a
^{-2r}c^{2\a br} ; K n^r ) ,$
%
\setcounter{equation}{9}
\begin{eqnarray}
\label{LP6} \max_{i\in\{0;1\}} \p_i \bigl(|\hat
x_n/x_{F_i,n} -1| \a_{{F_i}}^{2(1-r)}
c_{{F_i}}^r / w_{{F_i}} t_{i,n}^\star
\ge un^{-r} /2b \bigr) &\ge& \bigl(1-u^{1/r} /8n
\bigr)^{2n} /4 ,
\\
\label{LP7} \max_{i\in\{0;1\}} \p_i
\bigl(|x_{F_i,n} /\hat x_n -1| \a_{{F_i}}^{2(1-r)}
c_{{F_i}}^r / w_{{F_i}} t_{i,n}^\star
\ge un^{-r} /2b \bigr) &\ge& \bigl(1-u^{1/r} /8n
\bigr)^{2n} /4 ,
\end{eqnarray}
where $ \max_{i\in\{0;1\}} |t^\star_{i,n}-1|\to0 $ as $ n\to\infty$.
\end{theorem}

\section{Proofs}

Our approach to establishing lower bounds requires constructing two
distribution functions $ F_0 $ and $ F_1 ,$ where $ F_0 $ is a Pareto
d.f. and $ F_1\equiv F_{1,n} $ is a ``disturbed'' version of $ F_0 $.
We then apply Lemma \ref{LbH} that provides a non-asymptotic lower
bound to the accuracy of estimation when choosing between two close
alternatives.

The problem of estimating the tail index, the tail constant and $
x_{{F,n}} $ from $ X_1,\ldots,X_n $ is equivalent to the problem of estimating
$ \a_{F} ,$ $ c_{F} $ and quantiles from a sample $ Y_1,\ldots,Y_n $ of
i.i.d. positive r.v.s with the distribution
%
\begin{equation}
\label{T2} F(y) \equiv\p(Y\le y) = y^\a\ell(y)\qquad (y>0),
\end{equation}
where function $ \ell$ slowly varies at the origin.

We denote by $ {\cal F} $ the class of distributions obeying (\ref{T2}).
Note that $ \L(Y)\in{\cal F} $ if and only if $ \L(1/Y)\in{\cal H}.$
Obviously, a tail index estimator $ \a_n(X_1,\ldots,X_n) $ can be
considered an estimator $ \a_n(1/Y_1, \ldots, 1/Y_n) $ of index $ \a$
from the sample $ Y_1=1/X_1,\ldots,Y_n=1/X_n$, and vice versa.
The tradition of dealing with this equivalent problem stems from
\cite{H82}. We proceed with this equivalent formulation.

A counterpart to $ {\cal H}(b) $ is the following non-parametric class
of d.f.s on $ (0;\infty)$:
%
\begin{equation}
\label{T20} {\cal F}(b) = \Bigl\{ F\in{\cal F}\dvt \sup_{y<K^*(F)}
\bigl|c_{F}^{-1}y^{-\a_{F}}F(y)-1\bigr| y^{-b\a_{F}} < \infty
\Bigr\} ,
\end{equation}
where $ b>0 $ and $ K^*(F) $ is the right end-point of $ F$.
A d.f. $ F\in{\cal F}(b) $ obeys
\[
F(y) = c_{F} y^{\a_{F}} \bigl(1+\mathrm{O}\bigl(y^{b\a_{F}}\bigr)
\bigr)\qquad (y\to0),
\]
where $ \a_{F} = \lim_{y\downarrow0} (\ln F(y))/\ln y $ and
$ c_{F} = \lim_{y\downarrow0} y^{-\a_F} F(y).$

\begin{pf*}{Proof of Theorem \ref{LP-1}} Let $ h\in(0;c)$, and denote
\[
\a_0=\a, \qquad\a_1 = \a+\gamma,\qquad \gamma= h^{\a b} .
\]
We will employ the distribution functions $ F_0 $ and $ F_1,$ where
\begin{eqnarray*}
F_0(y) &=& (y/c)^\a\1\{0<y\le c\} ,
\\
F_1(y) &=& (h/c)^{-\gamma}(y/c)^{\a_1}\1\{0<y\le h
\}+(y/c)^\a\1\{ h<y\le c\} .
\end{eqnarray*}
The counterparts to these distributions are
\begin{eqnarray*}
\p_0(X > x) &=& (cx)^{-\a} \1\{x\ge1/c\} ,
\\
\p_1(X > x) &=& (cx)^{-\a} \1\{1/c \le x < 1/h\} +
c^{-\a}h^{-\gamma
}x^{-\a_1} \1\{x \ge1/h\}.
\end{eqnarray*}
It is easy to see that $ F_1\le F_0 $ and
%
\begin{equation}
\label{LP2} \a_{{F_0}}=\a, \qquad\a_{{F_1}}=\a_1 ,\qquad
c_{{F_0}}=c^{-\a} , \qquad c_{{F_1}}=c^{-\a}
h^{-\gamma}.
\end{equation}
Obviously, $ F_0\in{\cal F}(b).$ We now check that $ F_1\in{\cal F}(b).$

Since
\[
c_{F_1}^{-1}y^{-\a_1}F_1(y) =
y^{-\gamma}h^\gamma\qquad(h<y\le c),
\]
we have
%
\begin{equation}
\label{LP3} \sup_{0<y\le c} \bigl|1-c_{F_1}^{-1}y^{-\a_1}F_1(y)\bigr|y^{-b\a_1}
= \sup_{h<y\leq c} \bigl(1 - y^{-\gamma}h^\gamma
\bigr)y^{-b\a_1} .
\end{equation}
The right-hand side of (\ref{LP3}) takes on its maximum at $ y_0 = h
(1+\gamma/b\a_1)^{1/\gamma}$; the supremum is bounded by $
A:=\mathrm{e}^{1/\mathrm{e}\a
}/b\a$.
Note that $ \{F_0,F_1\}\subset{\cal D}_{b,A} .$

Let $ d^2_{ _H}(P_0;P_1) $ denote the Hellinger distance. It is easy to
check that
%
\begin{equation}
\label{LP5} d_{ _H}^2(F_0;F_1)
\le\gamma^{1/r} /8\a^2c^\a.
\end{equation}
According to Lemma \ref{LbH} below,
%
\begin{equation}
\label{LP10} \max_{i\in\{0;1\}} \p_i\bigl(|\hat
\a_n-\a_{F_i}|\ge\gamma/2\bigr) \ge \bigl(1-\gamma^{1/r}
/8\a^2c^\a\bigr)^{2n} /4 .
\end{equation}
Let $ \gamma=\gamma_n ,$ where
\[
\gamma_n \equiv\gamma_n(\a,c,v) = v \bigl(
\a^2c^\a/n\bigr)^r .
\]
Note that $ h<c $ as $ n > \a^2c^{-2b\a}v^{1/r}.$ From (\ref{LP10}),
%
\begin{equation}
\label{LP11} \max_{i\in\{0;1\}} \p_i \bigl(|\hat
\a_n/\a_{{F_i}}-1| \a_{{F_i}}^{r/b}
c_{{F_i}}^r n^r \ge vt_{n,i}/2 \bigr)
\ge\bigl(1-v^{1/r}/8n\bigr)^{2n} /4 ,
\end{equation}
where $ t_{n,0}=1 $ and $ t_{n,1}=1/f(\gamma),$ $ f(\gamma) =
(1+\gamma
/\a)^2\gamma^{\gamma/\a b} .$ Note that $ f(\gamma)\le1 $ as $
\gamma
\le \mathrm{e}^{-1-2b} $. Hence, $ t_{n,1}\ge1 $ as $ v\in[0;K n^r]$ and (\ref
{LP1}) follows.

Let $ \hat a_n $ be an arbitrary estimator of index $ a=1/\a$. Denote
$ a=a_0$. Since $ |a_0-a_1| = \gamma aa_1,$ Lemma \ref{LbH} yields
\[
\max_{i\in\{0;1\}} \p_i \bigl(|\hat a_n-a_{{F_i}}|
\ge\gamma aa_1/2 \bigr) \ge \bigl(1-\gamma^{1/r} /8
\a^2c^\a \bigr)^{2n} /4 .
\]
With $ \gamma= \gamma_n ,$ the left-hand side of this inequality is
\[
\max_{i\in\{0;1\}} \p_i \bigl(|\hat
a_n - a_{{F_i}}| \ge vn^{-r}a^{1-2r}a_1/2c_{{F_0}}^r
\bigr) = \max_{i\in\{0;1\}} \p_i \bigl(|\hat
a_n/a_{{F_i}} - 1|a_{{F_i}}^{-r/b}c_{{F_i}}^r
n^r \ge vt_{n,i}^+/2 \bigr) ,
\]
where $ t_{n,0}^+=(1+\gamma a)^r\ge1 $ and $ t_{n,1}^+=(1+\gamma
a)^{r/b} \gamma^{-r\gamma a/b}\ge1$, leading to (\ref{LP}).
\end{pf*}

\begin{pf*}{Proof of Corollary \ref{LP-C}} Note that
%
\begin{equation}
\label{LP9} \E\xi= \int_0^\infty\p(\xi\ge x)\,\mathrm{d}x
\end{equation}
for any non-negative r.v. $ \xi.$ Since
\[
\int_0^{z_n} \bigl(1-v^{1/r}/8n
\bigr)^{2n} \,\mathrm{d}v = 4^rr\Gamma(r)+\mathrm{o}(1)\qquad (n\to\infty)
\]
as $ z_n\to\infty,$ $ z_n=\mathrm{o}(n^r),$ (\ref{LP1}) and (\ref{LP9})
entail (\ref{LPC}).
\end{pf*}

\begin{pf*}{Proof of Theorem \ref{LP-2}} With $ F_0 $ and $ F_1 $
defined as above,
we have
\[
c_{F_1} - c_{F_0} = c^{-\a}\bigl(
\gamma^{-\gamma/\a b}-1\bigr) \ge c^{-\a} \gamma|\ln\gamma|/\a b .
\]
Using this inequality, (\ref{LP10}) and Lemma \ref{LbH}, we derive
\[
\max_{i\in\{0;1\}} \p_i \bigl(|\hat c_n-c_{{F_i}}|
\ge c^{-\a} \gamma|\ln\gamma|/2\a b \bigr) 
\ge \bigl(1-
\gamma^{1/r}/8\a^2c^\a \bigr)^{2n} /4 .
\]
Let $ \gamma\equiv\gamma(n) = (v\a^2c^\a/n)^r$. Then
\[
\max_{i\in\{0;1\}} \p_i \bigl(|\hat c_n-c_{{F_i}}|
\ge c_{{F_0}} \bigl(v\a^2c^\a/n
\bigr)^r r \ln(n/\ln n)/2\a b \bigr) \ge (1-v/8n )^{2n} /4 .
\]
Note that $ \a^2c^\a/\a^2_1c_{{F_1}}^{-1}
\ge t_n $ as $ v\in[0;\a^{-2}c^{-\a}\ln n]$. The result
follows.
\end{pf*}

\begin{pf*}{Proof of Theorem \ref{LP-3}} Denote
\[
x_i \equiv x_{F_i,n} ,\qquad y_i=1/x_i
.\vadjust{\goodbreak}
\]
Obviously, $ y_i $ is the quantile of $ \L_i(1/X)$. We find convenient
dealing with the equivalent problem of estimating quantiles of the
distribution of a random variable $ Y=1/X$.

With functions $ F_0, F_1 $ defined as above, it is easy to see that
%
\begin{equation}
\label{LP8} y_0 = cq_n^{1/\a} = c\kappa h ,\qquad
y_1 = c^{\a/\a_1} q_n^{1/\a_1}
h^{\gamma/\a_1} = y_0 (c\kappa)^{-\gamma/\a_1} ,
\end{equation}
where we put $ \kappa= q_n^{1/\a}/h .$ Note that $ y_1 = h(c\kappa
)^{1-\gamma/\a_1} .$ Hence $ y_i<h $ if $ c\kappa<1 $ ($i\in\{0,1\}$).

Denote
%
\begin{equation}
\label{gama} 
\gamma\equiv\gamma_n(\a,b,c) = u\bigl(
\a^2c^\a/n\bigr)^r .
\end{equation}
Then $ \kappa= s^{1/\a} (\a^2c^\a)^{-r/\a b}u^{-1/\a b} $ and
%
\begin{equation}
\label{LP13} c\kappa= u^{-1/\a b} s^{1/\a} c^{2r}
\a^{-2r/\a b} < 1
\end{equation}
by the assumption.

Using the facts that $ \mathrm{e}^{x}-1 \ge x\mathrm{e}^{x/2} $ and $ 1-\mathrm{e}^{-x} \ge
x\mathrm{e}^{-x/2} $ as $ x\ge0$, we derive
\begin{eqnarray*}
y_1 - y_0 &=& y_0 \bigl( (c
\kappa)^{-\gamma/\a_1} - 1 \bigr)
\\
&\ge& \gamma^{1+1/\a b} (c\kappa)^{1-\gamma/2\a_1} |\ln c\kappa |/
\a_1 .
\end{eqnarray*}
Hence, $ (y_1 - y_0)/y_0 \ge\gamma|\ln c\kappa|/\a_1 $ and
$ (y_1 - y_0)/y_1 = 1-(c\kappa)^{\gamma/\a_1}\ge\gamma|\ln c\kappa|
(c\kappa)^{\gamma/2\a_1}/\a_1$.
By Lemma \ref{LbH},
\[
\max_{i\in\{0;1\}} \p_i\bigl(|\hat y_n -
y_i| \ge \gamma^{1+1/\a b} (c\kappa)^{1-\gamma/2\a_1} |\ln c
\kappa|/2\a_1\bigr) \ge \bigl(1-\gamma^{1/r} /8
\a^2c^\a \bigr)^{2n} /4
\]
for any estimator $ \hat y_n $. Thus,
\[
\max_{i\in\{0;1\}} \p_i\bigl( |\hat y_n/y_i-1|
\ge\gamma\bigl|\ln(c\kappa )^{\a b}\bigr| t^\star_{n,i}/2b
\a_{{F_i}}^2 \bigr) \ge \bigl(1-\gamma^{1/r} /8\a
^2c^\a \bigr)^{2n} /4 ,
\]
where $ t^\star_{n,0} = 1/(1 + \gamma/\a)= 1/(1 + u\a_{{F_0}}^{-r/b}
c_{{F_0}}^{-r}n^{-r}) $ and $ t^\star_{n,1} =
(1 + \gamma/\a) (c\kappa)^{\gamma/2\a} .$ Taking into account
(\ref
{gama}) and (\ref{LP13}), we derive
\[
\max_{i\in\{0;1\}} \p_i \bigl( |\hat
y_n/y_i-1| \ge u \a_{{F_i}}^{2(r-1)}
c_{{F_i}}^{-r} n^{-r} \ln\bigl(u\a^{2r} /
s^bc^{2r\a b}\bigr) t^\star_{n,i}/2b \bigr)
\ge \bigl(1-u^{1/r} /8n \bigr)^{2n} /4 .
\]
This leads to (\ref{LP7}).

Recall that $ x_i=1/y_i$. From (\ref{LP8}),
\[
|x_1 - x_0| = |y_1 - y_0|/y_0y_1
\ge\gamma^{1-1/\a b} (c\kappa)^{-1+\gamma/2\a_1} |\ln c\kappa|/\a_1
.
\]
By Lemma \ref{LbH},
\[
\max_{i\in\{0;1\}} \p_i \bigl( |\hat
x_n-x_i| \ge\gamma^{1-1/\a b} (c
\kappa)^{-1+\gamma/2\a_1} |\ln c\kappa|/2\a_1 \bigr) \ge \bigl(1-
\gamma^{1/r}/8\a^2c^\a \bigr)^{2n} /4 .
\]
Hence,
\begingroup
\abovedisplayskip=7pt
\belowdisplayskip=7pt
\[
\max_{i\in\{0;1\}} \p_i \bigl( |\hat
x_n/x_i-1| \ge u\a_{{F_i}}^{2(r-1)}
c_{{F_i}}^{-r} n^{-r} \bigl|\ln\bigl(s^bc^{2r\a b}/
\a^{2r}u\bigr)\bigr|\tilde t_{n,i}/2b \bigr) \ge
\bigl(1-u^{1/r}/8n \bigr)^{2n} /4 ,
\]
%
where $ \tilde t_{n,0} = (c\kappa)^{\gamma/2\a}/(1 + \gamma/\a) =
(u^{-1/\a b} s^{1/\a} c^{2r} \a^{-2r/\a b})^{\gamma/2\a} / (1 + u\a
_{{F_0}}^{-r/b} c_{{F_0}}^{-r}n^{-r}) $ and $ \tilde t_{n,1}=1$.
The proof is complete.\vspace*{-2pt}
\end{pf*}

The next lemma presents a lower bound to the accuracy of choosing between
two ``close'' alternatives.

Let $ {\cal P} $ be an arbitrary class of distributions, and assume
that the quantity of interest, $ a_P,$ is an element of a metric space
$ ({\cal X},d)$.
An estimator $ \hat a $ of $ a_P $ is a measurable function of $
X_1,\ldots,X_n $ taking values in a subspace $ \{a_P\dvt P\in{\cal P}\} $ of
the metric space $ ({\cal X},d)$.

Examples of functionals $ a_P $ include
(a) $ a_{P_\theta}=\theta,$ where $ {\cal P} = \{P_\theta, \theta
\in\Theta\} $ is a parametric family of distributions ($\Theta
\subset
\R^m$);
(b) $ a_P=f_P,$ where $ f_P $ is the density of $ P $
with respect to a particular measure;
(c) $ a_P=P$.
A minimax lower bound over $ {\cal P} $ follows from a lower bound to
$ \max_{i\in\{0;1\}} \p_i(d(\hat a;a_{P_i}) \ge\delta),$ where $ P_0,
P_1\in{\cal P}.$\vspace*{-2pt}

\begin{lemma} \label{LbH}
Denote $ 2\delta= d(a_{P_1};a_{P_0})$. Then
%
\begin{equation}
\label{LP12} \max_{i\in\{0;1\}} \p_i\bigl(d(\hat
a;a_{P_i}) \ge\delta\bigr) \ge \bigl(1-d^2_{H}
\bigr)^{2n}/4 ,
\end{equation}
where $ d_{H} \equiv d_{H}(P_0;P_1) $ is the Hellinger
distance.\vspace*{-2pt}
\end{lemma}

There is considerable literature on techniques of deriving minimax
lower bounds of this kind (cf. \cite{Hu97,IH81,Tsy}). Classical results
include Fano's and Assuad's lemmas.
Inequality (\ref{LP12}) is sharper than Lemma 1 in \cite{Hu97}.
Another related result is Theorem 2.2 in \cite{Tsy}.\vspace*{-2pt}

\begin{pf*}{Proof of Lemma \ref{LbH}} Recall that
\[
d_{H}^2(P_0;P_1) = \frac12
\int \bigl(f_0^{1/2}-f_1^{1/2}
\bigr)^2 = 1 - \int\sqrt{f_0f_1} ,
\]
where $ f_i $ is a density of $ P_i $ with respect to a certain measure
(e.g., $ P_0+P_1$).

Let $ f_{i,n} $ denote the density of $ \L_i(X_1,\ldots,X_n),$ and
put $ a_i = a_{P_i}.$ By the triangle inequality, $ 2\delta\le
d(a_{P_0};\hat a) + d(\hat a;a_{P_1}).$ Therefore, $ 1\le\1_0+\1_1,$ where
\[
\1_0 = \1\bigl\{d(a_0;\hat a) \ge\delta\bigr\} ,\qquad
\1_1 = \1\bigl\{d(\hat a;a_1)\ge \delta\bigr\} .
\]
Using the definition of $ d_{ _H} $ and the Bunyakovskiy--Cauchy--Schwarz
inequality, we derive
\begin{eqnarray*}
\bigl(1-d^2_{ _H}\bigr)^n &=& \int
\sqrt{f_{0,n}f_{1,n}}
\\[-2pt]
&\le& \int\sqrt{f_{0,n}f_{1,n}} \1_0 + \int
\sqrt{f_{0,n}f_{1,n}} \1_1
\\[-2pt]
&\le& \p_0^{1/2}\bigl(d(\hat a;a_0) \ge\delta
\bigr) + \p_1^{1/2}\bigl(d(\hat a;a_1) \ge
\delta\bigr) .
\end{eqnarray*}
Hence $ (1-d^2_{ _H})^{2n} \le2 ( \p_0(d(\hat a;a_0) \ge\delta) +
\p_1(d(\hat a;a_1) \ge\delta)  )$, leading to (\ref{LP12}).\vadjust{\goodbreak}
\end{pf*}\endgroup


\section*{Acknowledgements}
The author is grateful to the Editor, the Associate Editor and two
referees for many helpful remarks.
Supported by a grant from the London Mathematical Society.



%

\printhistory

\end{document}